# Anomaly Detection Under Uncertainty Using Distributionally Robust Optimization Approach

Amir Hossein Noormohammadi[a,*], Seyed Ali MirHassani[a], Farnaz Hooshmand Khaligh[a]

[a] Department of Mathematics and Computer Science, Amirkabir University of Technology (Tehran Polytechnic), Tehran, Iran

## Abstract

Anomaly detection is defined as the problem of finding data points that do not follow the patterns of the majority. Among the various proposed methods for solving this problem, classification-based methods, including one-class Support Vector Machines (SVM) are considered effective and state-of-the-art. The one-class SVM method aims to find a decision boundary to distinguish between normal data points and anomalies using only the normal data. On the other hand, most real-world problems involve some degree of uncertainty, where the true probability distribution of each data point is unknown, and estimating it is often difficult and costly. Assuming partial distribution information such as the first and second-order moments is known, a distributionally robust chance-constrained model is proposed in which the probability of misclassification is low. By utilizing a mapping function to a higher dimensional space, the proposed model will be capable of classifying origin-inseparable datasets. Also, by adopting the kernel idea, the need for explicitly knowing the mapping is eliminated, computations can be performed in the input space, and computational complexity is reduced. Computational results validate the robustness of the proposed model under different probability distributions and also the superiority of the proposed model compared to the standard one-class SVM in terms of various evaluation metrics.

# 1 Introduction

The history of anomaly detection dates back to studies conducted by the statistics community at the beginning of the nineteen century [1]. Due to the importance of anomaly detection during the time, many researchers from various domains have studied this problem and a wide range of approaches have been introduced [2]. In [2], an anomalous data point is defined as one that does not follow the patterns of the majority of the data points considered normal data, and the anomaly detection problem is defined as finding such data points. Although anomalies and outliers are considered the same concept, in this research, outliers refer to those in the normal class in the training data set that deviate compared to other normal data points, due to various reasons such as instrument error, measurement error, and unrepresentative sampling. Figure 1 shows a two-dimensional data set containing anomalous data points. The data point sets $N_1$ and $N_2$ are normal, while the data points $A_1$ and $A_2$ and data point set $A_3$ are anomalous.

[*] Corresponding author.
E-mail addresses: a.h.noomohammadi@aut.ac.ir (A. H. Noormohammadi), a_mirhassani@aut.ac.ir (S. A. MirHassani), f.hooshmand.khaligh@aut.ac.ir (F. Hooshmand Khaligh).

Figure 1: Two-dimensional dataset with normal data $N_1$ and $N_2$ and anomalous data $A_1$, $A_2$, and $A_3$. Adopted from *[2]*

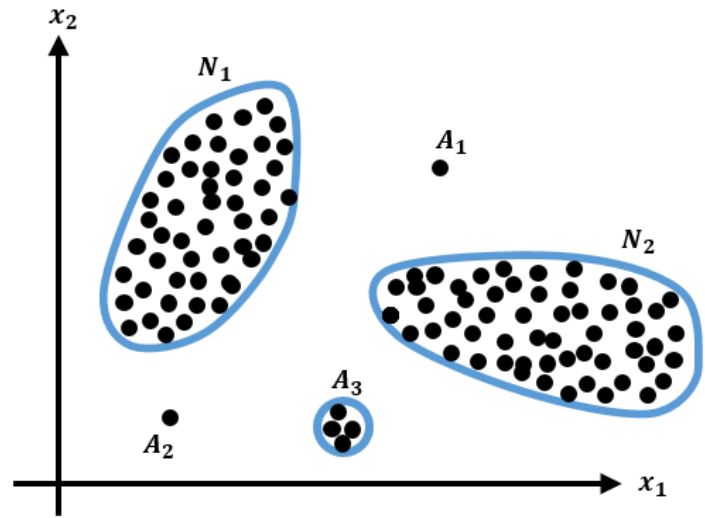

Anomaly detection is important because the presence of anomalous data points can significantly degrade the performance of machine learning algorithms. When anomalies exist in the training data, they distort the learning process, resulting in models that do not generalize well on normal data.

Anomaly detection encompasses a wide range of applications including healthcare monitoring to identify irregularities indicative of emerging health issues or atypical patient conditions. fault detection to understand abnormal behavior in an industrial production process [3]; fraud detection to identify criminal activities in financial organizations [4]; intrusion detection to monitor network traffic and detect attacks [5] and [6].

A vast range of methods have been proposed to address the anomaly detection problem, which can be categorized into 5 groups; statistical methods including parametric methods like Gaussian distribution-based [7], regression-based [8], distribution mixture-based [9], and nonparametric methods like kernel density estimation-based [10]; classification-based methods including autoencoder neural networks [11], one-class SVM [12], support vector data description [13], isolation forest [14]; density-based methods such as local outlier factor [15]; distance-based methods like k-nearest neighbors; and clustering-based methods including K-means [16] and DBSCAN [17].

In real-world applications, due to the high cost and in some cases infeasibility of collecting anomalies, very few anomalous data points are available. For example, in manufacturing quality monitoring, a fault detector determines if an item is defective or not by examining its features and comparing it to the majority of manufactured items. Clearly, producing faulty items is very costly [13]. Also, due to anomaly diversity, the observed anomalies may not represent future anomalies, hence the trained model may fail to detect new types of anomalies. Therefore, a method is needed that can perform anomaly detection using only the normal data. The one-class SVM method allows training a model using only the normal data points. In addition, due to the extensive use of SVM-based methods and their intuitive decision boundary description, the one-class SVM is considered as the basic method for solving the anomaly detection problem in this research.

The one-class SVM model was first introduced by Schölkopf et al. [12] in 2001 with the goal of finding a hyperplane to separate the normal data and anomalies by maximizing the distance of the hyperplane to the origin. Many researchers have attempted to make improvements to the basic model. In [18] by using the Ramp loss function instead of the Hinge loss function, the impact of outlier data points in the training data on the decision boundary is decreased. Also by reducing the number of support vectors, the obtained solution becomes sparser and the model training time decreases [17]. Since the one-class SVM performance heavily relies on the hyperparameter values, proper hyperparameter tuning is very important. In [18], an edge pattern recognition approach has been used to tune the penalty coefficient and Gaussian kernel width hyperparameters.

In all the aforementioned methods, it is assumed that the data points are known exactly. However, in practice, due to sampling, modeling, and measurement errors, the provided feature values are approximations of the actual values or have deviations from the actual values due to noise. Also, in some cases, some feature values are missed and replaced with other values. As a result, the feature vectors encompass some degree of uncertainty, making the study of real-world classification data sets more challenging [19]. In addition to the uncertainty in feature vectors, in practice, the true probability distribution of a data point is also unknown, and estimating it is often difficult and costly [19].

To the best of our knowledge, the uncertainty in feature vectors and their data distributions has not been previously studied in the anomaly detection problem, and uncertainty has only been examined in special cases – such as additive bounded noise – using robust optimization, which is conservative due to ignoring the hidden distribution information embedded in the data sets [20] and [21]. Also, most studies on handling uncertainty have only focused on the linear case of the one-class SVM, where the data set is origin-separable, while in most practical applications, the data set is origin-inseparable separable.

Motivated by the research in [18] on incorporating uncertainty into support vector machine modeling, this paper adopts a distributionally robust optimization approach to address uncertainty, with the following innovations: Uncertainty is simultaneously considered in the data points and their probability distributions in the one-class SVM model. A nonlinear model is provided to classify origin-inseparable datasets. The kernel idea is utilized to reduce computational complexity. The proposed model has been tested on several simulated datasets, and its robustness to uncertainty in different distributions and its superiority over the standard one-class SVM in terms of various evaluation metrics have been shown.

The organization of this paper is as follows: Section 2 reviews the standard one-class SVM model. Section 3 describes the problem and formulates the one-class SVM model under uncertainty, reformulates it as a deterministic model, and presents a kernel-based model. Section 4 contains the evaluation metrics and computational results. Section 5 concludes this paper.

## 2 One-Class SVM
In this section, the standard one-class SVM model is reviewed, assuming the dataset is origin-inseparable and the feature vector of each data point is known exactly.

### 2.1 Problem Description
A dataset $D_0$ is called origin-separable if there exists a vector w such that $\mathbf{w}^T\mathbf{x}^{(i)} > 0 \ \ \forall i = 1, \dots, l$. If $D_0$ is origin-separable, then there exists a unique hyperplane $f_{\mathbf{w},b}(\mathbf{x}) = \{\mathbf{x} \in \mathbb{R}^n : \mathbf{w}^T\mathbf{x} - b = 0\}$ with $b > 0$ that separates all points from the origin with the maximum distance to the origin among all separating hyperplanes [10]. For simplicity, the hyperplane is denoted by $\mathbf{w}^T\mathbf{x} - b = 1$. In most practical applications, the dataset is origin-inseparable and needs to be mapped to a new space such that it becomes origin-separable. If $D_0$ is origin-inseparable, a mapping function $\boldsymbol{\phi}: \mathbb{R}^n \to \mathbb{R}^d$ is defined to map the feature vectors from the input space $\mathbb{R}^n$ to a higher dimensional called feature space $\mathbb{R}^d$ where $d > n$, such that the dataset becomes origin-separable in the feature space.

Suppose an origin-inseparable unlabeled training dataset $\mathcal{D}_0 = \{\mathbf{x}^{(i)}: i = 1, \ldots, l\}$, including normal data points with some outliers, where $\mathbf{x}^{(i)} = \left[\mathbf{x}_1^{(i)}, \ldots, \mathbf{x}_n^{(i)}\right]^T$ represents the feature vector of the $i^{th}$ data point in n-dimensional real space $\mathbb{R}^n$. The goal is to find a hyperplane $f_{\mathbf{w},b}(\mathbf{x}) = \mathbf{w}^T\mathbf{x} - b = 0$ such that most data points lie on the side opposite to the origin with maximum distance to the origin. The distance of the hyperplane to the origin is called the margin. Since the inner product of the normal data points and the origin is zero, the normal data has the least similarity with the origin, hence maximizing the margin results in the minimum generalization error (test data error). The vector $\mathbf{w}$ is called the normal vector and the scalar $b$ is called the bias.

## 2.2 Nonlinear One-Class SVM Formulation

The primal model of nonlinear soft margin one-class SVM is formulated as follows:

$$\min_{\mathbf{w},b\geq 0,\boldsymbol{\xi}} \frac{1}{2}\|\mathbf{w}\|_2^2 - b + \frac{1}{l\nu}\sum_{i=1}^{l}\xi_i$$

$$\mathbf{w}^T\boldsymbol{\phi}(\mathbf{x}^{(i)}) \geq b - \xi_i \quad i = 1, \ldots, l \qquad (1)$$

$$\xi_i \geq 0 \quad \forall i \in \{1, \ldots, l\}$$

Figure 2 shows a geometric illustration of the hyperplane obtained from the one-class SVM model on an origin-separable dataset.

Figure 2: Geometric intuition of the decision boundary obtained from one-class SVM on an origin-separable dataset. Adopted from [18]

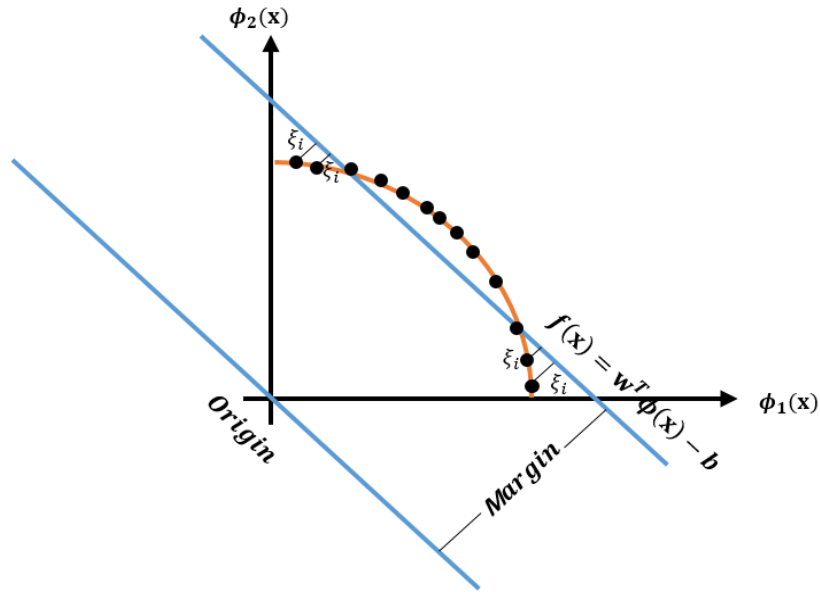

In the primal model (1), the mapping function $\boldsymbol{\phi}$ is explicitly utilized, however, there is no need to know $\boldsymbol{\phi}$ in practice. Hence, we derive the dual formulation of model (1). By introducing the Lagrange multipliers $\alpha \geq 0$ and $\beta \geq 0$, the Lagrangian is formed as:

$$\mathcal{L}(\mathbf{w}, \boldsymbol{\xi}, b, \alpha, \beta) = \frac{1}{2}\|\mathbf{w}\|_2^2 - b + \frac{1}{l\nu}\sum_{i=1}^{l}\xi_i - \sum_{i=1}^{l}\alpha_i\left(\mathbf{w}^T\boldsymbol{\phi}(\mathbf{x}^{(i)}) - b + \xi_i\right) - \sum_{i=1}^{l}\beta_i\xi_i \qquad (2)$$

According to the first-order necessary optimality condition for convex functions, the optimal solution of $\min_{\mathbf{w}\in\mathbb{F}, \boldsymbol{\xi}\geq 0, \rho} \mathcal{L}(\mathbf{w}, \boldsymbol{\xi}, b, \boldsymbol{\alpha}, \boldsymbol{\beta})$ is obtained by taking the partial derivatives of the Lagrangian (2) with respect to the variables $\mathbf{w}$, $\boldsymbol{\xi}$, and $b$ and solving the following system:

$$\mathbf{w} = \sum_{i=1}^{l} \alpha_i \boldsymbol{\phi}(\mathbf{x}^{(i)}) \tag{3}$$

$$\sum_{i=1}^{l} \alpha_i = 1 \tag{4}$$

$$\alpha_i = \frac{1}{l\nu} - \beta_i \quad i = 1, \dots, l \tag{5}$$

By substituting equations (3)-(5) into the Lagrangian equation (2), the Lagrange dual function becomes:

$$\mathcal{G}(\boldsymbol{\alpha}, \boldsymbol{\beta}) = \frac{1}{2}\sum_{i=1}^{l}\sum_{i'=1}^{l} \alpha_i \alpha_{i'} \boldsymbol{\phi}(\mathbf{x}^{(i)})^{\mathrm{T}} \boldsymbol{\phi}(\mathbf{x}^{(i')}) - b + \frac{1}{\nu l}\sum_{i=1}^{l}\xi_i - \sum_{i=1}^{l}\sum_{i'=1}^{l}\alpha_i \alpha_{i'}\boldsymbol{\phi}(\mathbf{x}^{(i)})^{\mathrm{T}}\boldsymbol{\phi}(\mathbf{x}^{(i')})$$
$$+ b\sum_{i=1}^{l}\alpha_i - \sum_{i=1}^{l}\alpha_i \xi_i - \sum_{i=1}^{l}\beta_i \xi_i$$

Finally, the dual problem of model (1) is formulated as:

$$\min_{\alpha \geq 0} \frac{1}{2}\sum_{i=1}^{l}\sum_{i'=1}^{l} \alpha_i \alpha_{i'} \boldsymbol{\phi}(\mathbf{x}^{(i)})^{\mathrm{T}} \boldsymbol{\phi}(\mathbf{x}^{(i')})$$

$$\sum_{i=1}^{l} \alpha_i = 1 \tag{6}$$

$$0 \leq \alpha_i \leq \left(\frac{1}{\nu l}\right) \quad i = 1, \dots, l$$

## 2.3 Formulation of Kernel-based One-Class SVM

Due to the increase in the number of variables in model (6) from using the mapping $\boldsymbol{\phi}$ and the computational complexity of the inner product $\boldsymbol{\phi}(\mathbf{x}^{(i)})^{\mathrm{T}}\boldsymbol{\phi}(\mathbf{x}^{(i')})$, especially in an infinite dimensional feature, the kernel idea is used to address these issues. More precisely, the kernel function replaces the inner product in the feature space.

The kernel function $k: \mathbb{R}^n \times \mathbb{R}^n \to \mathbb{R}$ is a similarity function that measures the similarity between the vectors x and x′. Some of the well-known kernels include linear kernel, polynomial kernel, and Gaussian kernel. In this paper, we assume the kernel is valid, meaning there exists a mapping such as $\boldsymbol{\phi}: \mathbb{R}^n \to \mathbb{R}^d$ such that $k(\mathbf{x}, \mathbf{x}') = \boldsymbol{\phi}(\mathbf{x})^T \boldsymbol{\phi}(\mathbf{x}')$. Linear, polynomial, and

Gaussian kernels are valid kernels. In this paper, the valid Gaussian kernel is utilized which has the following equation:

$$k\left(\mathbf{x}^{(i)}, \mathbf{x}^{(i')}\right) = e^{-\gamma \left\|\mathbf{x}^{(i)} - \mathbf{x}^{(i')}\right\|_2^2} \tag{7}$$

The Gaussian kernel (7) has only one parameter $\gamma$; The corresponding feature space has infinite dimension, the mapped feature vectors in this space have acute angles, and all these vectors lie on a hypersphere with unit radius centered at the origin, therefore, as mentioned in [10], the Gaussian kernel guarantees that the dataset becomes origin-inseparable in the feature space.

As can be seen from equation (7), the kernel $k\left(\mathbf{x}^{(i)}, \mathbf{x}^{(i')}\right)$ and the corresponding inner product $\boldsymbol{\phi}\left(\mathbf{x}^{(i)}\right)^T \boldsymbol{\phi}\left(\mathbf{x}^{(i')}\right)$ approach 1 when $\mathbf{x}^{(i)}$ and $\mathbf{x}^{(i')}$ are close, indicating high similarity between the mapped vectors, while $k\left(\mathbf{x}^{(i)}, \mathbf{x}^{(i')}\right)$ approaches 0 when $\mathbf{x}^{(i)}$ and $\mathbf{x}^{(i')}$ diverge, indicating dissimilarity between the mapped vectors useful for differentiating normal and anomalous points in the feature space.

Now instead of the inner product $\boldsymbol{\phi}\left(\mathbf{x}^{(i)}\right)^T \boldsymbol{\phi}\left(\mathbf{x}^{(i')}\right)$ in model (7), $k\left(\mathbf{x}^{(i)}, \mathbf{x}^{(i')}\right)$ can be used, therefore the Kernel-based one-class SVM is expressed as the following quadratic optimization problem:

$$\begin{aligned}
&\min_{\alpha \geq 0} \frac{1}{2} \sum_{i=1}^{l} \sum_{i'=1}^{l} \alpha_i \alpha_{i'} k\left(\mathbf{x}^{(i)}, \mathbf{x}^{(i')}\right) \\
&\sum_{i=1}^{l} \alpha_i = 1 \\
&0 \leq \alpha_i \leq \left(\frac{1}{\nu l}\right) \quad i = 1, \dots, l
\end{aligned} \tag{8}$$

By solving model (8) and finding the solution $\alpha^*$, according to equation (2) the hyperplane coefficients $\mathbf{w}^*$ are obtained as:

$$\mathbf{w}^* = \sum_{i=1}^{l} \alpha_i^* \boldsymbol{\phi}\left(\mathbf{x}^{(i)}\right) \tag{9}$$

For the data point $\mathbf{x}_{sv}$ that its corresponding Lagrange multiplier $\alpha_{sv}^*$ satisfies $0 < \alpha_{sv}^* < \frac{1}{l\nu}$, the optimal value of $b$ is obtained as:

$$b^* = \sum_{i=1}^{l} \alpha_i^* k\left(\mathbf{x}^{(i)}, \mathbf{x}_{sv}\right) \tag{10}$$

And the decision function in terms of the dual variables is:

$$f(\mathbf{x}) = \sum_{i=1}^{l} \alpha_i^* k(\mathbf{x}^{(i)}, \mathbf{x}) - b^* \tag{11}$$

As proven in [12], since the Gaussian kernel (7) is used in model (8), the dataset is origin-inseparable, and according to the support vector hyperplane theorem [22], model (8) always has a unique solution. Any data point that $\mathbf{x}^{(i)}$ that its corresponding Lagrange multiplier $\alpha_i^*$ is positive is called a support vector. According to the KKT conditions, only these points contribute in determining the hyperplane coefficients.

As can be seen, the number of variables in model (8) is equal to the number of training data points, therefore in problems with a high number of features such as when the dataset is origin-inseparable and needs to be mapped in a higher dimensional space, solving model (9) can have lower complexity.

## 3. One-Class SVM Under Uncertainty

This section aims to explore the efficacy of distributionally robust chance-constrained one-class SVM with uncertain input data and probability distributions specified by the first- and second-order moments information. To address this uncertainty, distributionally robust constraints or ambiguous chance constraints have been developed and adopted to represent a conservative approximation of the original problem (8). The approach used in this paper ensures that the model can provide a classifier that has the best performance even in the worst case over a set of distributions. Also, considering the efficiency of the kernel idea in reducing computational complexity and eliminating the need for explicit use of feature mapping, a kernel-based reformulated model is presented.

### 3.1 Assumptions and Problem Statement

Assume an origin-inseparable unlabeled training dataset $\mathcal{D}_1 = \{\tilde{\mathbf{x}}^{(i)} \in \mathbb{R}^n : i = 1, \dots, m\}$ is given that contains the normal data points with some outliers, where $\tilde{\mathbf{x}}^{(i)} = [\tilde{x}_1^{(i)} \quad \dots \quad \tilde{x}_n^{(i)}]^T$ represents the random feature vector of the $i^{th}$ data point following the underlying probability distribution $F_U^{(i)}$. Let $F_U^{(i)}$ be mutually independent for $i = 1, \dots, l$, meaning for any $S \subseteq \{1, \dots, l\}$ we have $\mathbb{P}_{\cap_{i \in S} F_U^{(i)}} = \prod_{i \in S} \mathbb{P}_{F_U^{(i)}}$ where $\mathbb{P}_{F_U^{(i)}}$ is the probability measure corresponding to the probability distribution $F_U^{(i)}$. Assume that true probability distribution $F_U^{(i)}$ is unknown, but its two first moments, i.e., mean vector $\boldsymbol{\mu}^{(i)} \triangleq \mathbb{E}_{F_U^{(i)}}[\tilde{\mathbf{x}}^i]$ and covariance matrix $\Sigma^{(i)} \triangleq \mathbb{E}_{F_U^{(i)}}\left[(\tilde{\mathbf{x}}^i - \mathbb{E}_{F_i}[\tilde{\mathbf{x}}^i])(\tilde{\mathbf{x}}^i - \mathbb{E}_{F_i}[\tilde{\mathbf{x}}^i])^T\right]$ are known a priori. Since the true probability distribution is not known, the chance constraint corresponding to $i^{th}$ data is imposed over an uncertain set containing all distributions with mean $\boldsymbol{\mu}^{(i)}$ and covariance matrix $\Sigma^{(i)}$, called the ambiguity set of the $i^{th}$ data point, which is defined as follows:

$$\mathcal{P}_i\left(\tilde{\mathbf{x}}^{(i)}; \boldsymbol{\mu}^{(i)}, \Sigma^{(i)}\right) = \left\{ F^{(i)} \middle| \begin{array}{c} \mathbb{P}_{F^{(i)}}\left(\tilde{\mathbf{x}}^{(i)} \in \Xi^{(i)}\right) = 1 \\ \mathbb{E}_{F^{(i)}}\left[\tilde{\mathbf{x}}^{(i)}\right] = \boldsymbol{\mu}^{(i)} \\ \mathbb{E}_{F^{(i)}}\left[\left(\tilde{\mathbf{x}}^{(i)} - \boldsymbol{\mu}^{(i)}\right)\left(\tilde{\mathbf{x}}^{(i)} - \boldsymbol{\mu}^{(i)}\right)^T\right] = \Sigma^{(i)} \end{array} \right\} \tag{12}$$

Where $\Xi^{(i)}$ is the support set of the random variable $\tilde{\mathbf{x}}^{(i)}$, and $\mathbb{E}_{F^{(i)}}[,]$ is the expected value function of distribution $F_U^{(i)}$. The first constraint in the set (12) guarantees the property $\mathbb{P}_{F^{(i)}}\left(\Xi^{(i)}\right) =$

1 of the probability measure $\mathbb{P}_{F^{(i)}}$. The second and third constraints indicate that the first two moments of the probability distribution must equal to the given mean and covariance matrix of the $i^{th}$ data.

The goal is to find a hyperplane $f_{\mathbf{w},b}(\mathbf{x}) = \mathbf{w}^T \boldsymbol{\phi}(\mathbf{x}) - b = 0$ in the feature space with a maximum margin such that the given individual chance constraints for separating most normal data points from the origin are satisfied.

## 3.2 Nonlinear Distributionally Robust Individual Chance-Constrained One-Class SVM Formulation

The one-class SVM model under uncertainty is defined as:

$$\min_{\mathbf{w},b \geq 0, \boldsymbol{\xi}} \frac{1}{2}\|\mathbf{w}\|_2^2 - b + \frac{1}{l\nu}\sum_{i=1}^{l} \xi_i$$

$$\mathbf{w}^T \boldsymbol{\phi}(\tilde{\mathbf{x}}^{(i)}) \geq b - \xi_i \quad i = 1, \dots, l$$

$$\xi_i \geq 0 \quad \forall i \in \{1, \dots, l\}$$

(13)

One approach to dealing with uncertainty is to use individual chance constraints which guarantee an upper bound on the probability of misclassification of each uncertain training data point as a random vector. Assume an error level $\alpha$ is given; We aim to find a hyperplane that each uncertain training data point $\tilde{\mathbf{x}}^{(i)}$ has at most a probability $\alpha$ under the true distribution $F_U^{(i)}$ to individually violate the classification constraint, in other words, it gets misclassified. For this purpose, inspired by the proposed model [23], the one-class SVM model with individual chance constraints is defined as:

$$\min_{\mathbf{w},b \geq 0, \boldsymbol{\xi}} \frac{1}{2}\|\mathbf{w}\|_2^2 - b + \frac{1}{l\nu}\sum_{i=1}^{l} \xi_i \tag{14}$$

$$\mathbb{P}_{F_U^{(i)}}\left[\mathbf{w}^T \boldsymbol{\phi}(\tilde{\mathbf{x}}^{(i)}) \leq b - \xi_i\right] \leq \alpha \quad \forall i \in \{1, \dots, l\} \tag{15}$$

$$\xi_i \geq 0 \quad \forall i \in \{1, \dots, l\} \tag{16}$$

Where $0 < \alpha < 1$ indicates the error level and affects the conservativeness of the classification. The chance constraint (15) provides an upper bound on the probability of misclassification of $\tilde{\mathbf{x}}^{(i)}$.

As mentioned, in practice, the data distribution is not exactly known or difficult to obtain, and only partial information such as two first moments of the distribution is given. Therefore, instead of relying on a single estimate of the true distribution, considering the ambiguity set (13) consisting of distributions with the same two first moments as the given information, as well as considering the worst-case misclassification over such distributions, the nonlinear one-class SVM model with distributionally robust individual chance constraints inspired by [23] is defined as:

$$\min_{\mathbf{w},b \geq 0, \boldsymbol{\xi}} \frac{1}{2}\|\mathbf{w}\|_2^2 - b + \frac{1}{l\nu}\sum_{i=1}^{l} \xi_i \tag{17}$$

$$\sup_{F_i \in \mathcal{P}_i(\tilde{\mathbf{x}}^{(i)};\boldsymbol{\mu}^{(i)},\Sigma^{(i)})} \mathbb{P}_{F_U^{(i)}}\left[\mathbf{w}^T\boldsymbol{\phi}(\tilde{\mathbf{x}}^{(i)}) \leq b - \xi_i\right] \leq \alpha \quad \forall i \in \{1,\dots,l\} \tag{18}$$

$$\xi_i \geq 0 \quad \forall i \in \{1,\dots,l\} \tag{19}$$

It is observed that if the distributionally robust individual chance constraints (18) are satisfied, the individual chance constraints (15) are also satisfied.

As shown later, in the case of uncertainty, each training data point is not just a single point but takes its values from an ellipsoidal uncertainty set (the size of the ellipsoid depends on error level α). In this case, we aim to find a hyperplane so that not only each data point but also the ellipsoid containing that data point is classified as correctly as possible.

### 3.3 Reformulation of the Uncertain One-Class SVM into a Deterministic Model

Chance-constrained programs are generally nonconvex problems and solving such problems is challenging [19]. For this purpose, a reformulation of the model (17)-(19) into a deterministic second-order cone programming (SOCP) model is provided; multivariate Chebyshev inequality [24] is utilized for this reformulation [23].

Multivariate Chebyshev Inequality: Assume $\tilde{\mathbf{x}}$ is a random variable with mean $\boldsymbol{\mu}$ and covariance matrix $\Sigma$. For any closed convex set $S$, the probability that the random variable $\tilde{\mathbf{x}}$ takes values from the set $S$ has an upper bound given by:

$$\mathbb{P}_F\{\tilde{\mathbf{x}} \in S\} \leq \frac{1}{1+d^2} \quad \forall F \in \mathcal{P}(\tilde{\mathbf{x}};\boldsymbol{\mu},\Sigma)$$

Equivalently we have:

$$\sup_{F \in \mathcal{P}(\tilde{\mathbf{x}};\boldsymbol{\mu},\Sigma)} \mathbb{P}_F\{\tilde{\mathbf{x}} \in S\} = \frac{1}{1+d^2}$$

where $d^2 = \inf_{\tilde{\mathbf{x}} \in S}(\tilde{\mathbf{x}} - \boldsymbol{\mu})^T \Sigma^{-1}(\tilde{\mathbf{x}} - \boldsymbol{\mu})$.

Now assume the first and second-order moments of the mapped data $\boldsymbol{\phi}(\tilde{\mathbf{x}}^{(i)}) \in \mathbb{R}^d$ are exactly known as the mean vector $\boldsymbol{\mu}_{\boldsymbol{\phi}}^{(i)} \in \mathbb{R}^d$ and the covariance matrix $\Sigma_{\boldsymbol{\phi}}^{(i)} \succcurlyeq 0$, respectively. By considering the halfspace $S^{(i)} = \{\boldsymbol{\phi}(\tilde{\mathbf{x}}^{(i)}): \mathbf{w}^T\boldsymbol{\phi}(\tilde{\mathbf{x}}^{(i)}) \leq b - \xi_i\}$ as a closed convex set corresponding to the $i^{th}$ data, Multivariate Chebyshev inequality for the $i^{th}$ random vector becomes:

$$\sup_{F_i \in \mathcal{P}(\tilde{\mathbf{x}}^{(i)};\boldsymbol{\mu}^{(i)},\Sigma^{(i)})} \mathbb{P}_{F_i}\{\mathbf{w}^T\boldsymbol{\phi}(\tilde{\mathbf{x}}^{(i)}) \leq b - \xi_i\} = \frac{1}{1+d_i^2}$$

where $d_i^2 = \inf_{\mathbf{w}^T\boldsymbol{\phi}(\tilde{\mathbf{x}}^{(i)}) \leq b - \xi_i}\left(\boldsymbol{\phi}(\tilde{\mathbf{x}}^{(i)}) - \boldsymbol{\mu}_{\boldsymbol{\phi}}^{(i)}\right)^T \Sigma_{\boldsymbol{\phi}}^{(i)^{-1}}\left(\boldsymbol{\phi}(\tilde{\mathbf{x}}^{(i)}) - \boldsymbol{\mu}_{\boldsymbol{\phi}}^{(i)}\right)$.

If $\mathbf{w}^T\boldsymbol{\mu}_{\boldsymbol{\phi}}^{(i)} \leq b - \xi_i$, then we have:

$$d_i^2 = \inf_{\mathbf{w}^T\boldsymbol{\phi}(\tilde{\mathbf{x}}^{(i)}) \leq b - \xi_i}\left(\boldsymbol{\phi}(\tilde{\mathbf{x}}^{(i)}) - \boldsymbol{\mu}_{\boldsymbol{\phi}}^{(i)}\right)^T \Sigma_{\boldsymbol{\phi}}^{(i)^{-1}}\left(\boldsymbol{\phi}(\tilde{\mathbf{x}}^{(i)}) - \boldsymbol{\mu}_{\boldsymbol{\phi}}^{(i)}\right)$$

$$= \left(\boldsymbol{\mu}_{\boldsymbol{\phi}}^{(i)} - \boldsymbol{\mu}_{\boldsymbol{\phi}}^{(i)}\right)^T \left(\Sigma_{\boldsymbol{\phi}}^{(i)}\right)^{-1} \left(\boldsymbol{\mu}_{\boldsymbol{\phi}}^{(i)} - \boldsymbol{\mu}_{\boldsymbol{\phi}}^{(i)}\right) = 0$$

And hence $\sup_{F_i \in \mathcal{P}\left(\tilde{\mathbf{x}}^{(i)}; \boldsymbol{\mu}^i, \Sigma^i\right)} \mathbb{P}_{F_i}\left\{\mathbf{w}^T \tilde{\mathbf{x}}^{(i)} \leq b - \xi_i\right\} = 1$ but according to constraint (18) and upper bound $\alpha < 1$, this result is contradicted.

If $\mathbf{w}^T \boldsymbol{\mu}_{\boldsymbol{\phi}}^{(i)} > b - \xi_i$, let $\mathbf{u}_i = \left(\Sigma_{\boldsymbol{\phi}}^{(i)}\right)^{-\frac{1}{2}} \left(\tilde{\mathbf{x}}^{(i)} - \boldsymbol{\mu}_{\boldsymbol{\phi}}^{(i)}\right)$, $\mathbf{v}_i = \left(\Sigma_{\boldsymbol{\phi}}^{(i)}\right)^{\frac{1}{2}} \mathbf{w}$, and $\gamma_i = b - \xi_i - \mathbf{w}^T \boldsymbol{\mu}_{\boldsymbol{\phi}}^{(i)}$, then $d_i^2 = \inf_{\mathbf{v}_i^T \mathbf{u}_i \leq \gamma_i} \mathbf{u}_i^T \mathbf{u}_i$. It is a convex quadratic programming with linear constraints with respect to the vector variable $\mathbf{u}_i$; therefore, the KKT conditions are necessary and sufficient for optimality. By introducing the Lagrange multipliers $\lambda_i \geq 0$, the Lagrangian becomes:

$$\mathcal{L}(\mathbf{u}_i, \lambda_i) = \mathbf{u}_i^T \mathbf{u}_i + \lambda_i (\mathbf{v}_i^T \mathbf{u}_i - \gamma_i)$$

Now the KKT conditions become:

Stationarity Conditions:

$$\frac{\partial \mathcal{L}(\mathbf{u}_i, \lambda)}{\partial \mathbf{u}_i} = 2\mathbf{u}_i^* + \lambda_i^* \mathbf{v}_i = 0 \tag{20}$$

Primal Feasibility Conditions:

$$\mathbf{v}_i^T \mathbf{u}_i^* \leq \gamma_i \tag{21}$$

Dual Feasibility Conditions:

$$\lambda_i^* \geq 0 \tag{22}$$

Complementary Slackness Conditions:

$$\lambda_i^* (\mathbf{v}_i^T \mathbf{u}_i^* - \gamma_i) = 0 \tag{23}$$

From equation (20) it follows that $\mathbf{u}_i^* = -\frac{\lambda_i^*}{2} \mathbf{v}_i$, also, $\lambda_i^* > 0$ (if $\lambda_i^* = 0$, then $\mathbf{u}_i^* = 0$. According to equation (21), $\mathbf{v}_i^T \mathbf{u}_i^* \leq \gamma_i$ which would imply $\gamma_i \geq 0$, contradicting the assumption $\gamma_i < 0$). As a result, according to equation (23), we have $\mathbf{v}_i^T \mathbf{u}_i^* = \gamma_i$ which gives $(\mathbf{u}_i^*)^T \mathbf{u}_i^* = \frac{\gamma_i^2}{\mathbf{v}_i^T \mathbf{v}_i}$. Therefore

$$d_i^2 = (\mathbf{u}_i^*)^T \mathbf{u}_i^* = \frac{\gamma_i^2}{\mathbf{v}_i^T \mathbf{v}_i} = \frac{\left(\mathbf{w}^T \boldsymbol{\mu}_{\boldsymbol{\phi}}^{(i)} + \xi_i - b\right)^2}{\mathbf{w}^T \Sigma_{\boldsymbol{\phi}}^{(i)} \mathbf{w}} \tag{24}$$

Now according to multivariate Chebyshev inequality, constraints (18) and (24), and the inequality $\alpha > 0$ it follows that:

$$\frac{1-\alpha}{\alpha} \leq d_i^2 = \frac{\left(\mathbf{w}^T \boldsymbol{\mu}_{\boldsymbol{\phi}}^{(i)} + \xi_i - b\right)^2}{\mathbf{w}^T \Sigma_{\boldsymbol{\phi}}^{(i)} \mathbf{w}}.$$

Or equivalently:

$$\sqrt{\frac{1-\alpha}{\alpha}} \left\| \left(\Sigma_{\boldsymbol{\phi}}^{(i)}\right)^{\frac{1}{2}} \mathbf{w} \right\|_2 \leq \mathbf{w}^T \boldsymbol{\mu}_{\boldsymbol{\phi}}^{(i)} + \xi_i - b$$

This is a second-order cone constraint. According to the obtained constrain, the reformulation of model (17)-(19) into a deterministic second-order cone programming model is as follows:

$$\min_{\mathbf{w}, b \geq 0, \xi} \frac{1}{2} \|\mathbf{w}\|_2^2 - b + \frac{1}{l\nu} \sum_{i=1}^{l} \xi_i \tag{25}$$

$$\sqrt{\frac{1-\alpha}{\alpha}} \left\| \left(\Sigma_{\boldsymbol{\phi}}^{(i)}\right)^{\frac{1}{2}} \mathbf{w} \right\|_2 \leq \mathbf{w}^T \boldsymbol{\mu}_{\boldsymbol{\phi}}^{(i)} + \xi_i - b \quad i = 1, \ldots, l \tag{26}$$

$$\xi_i \geq 0 \quad \forall i \in \{1, \ldots, l\} \tag{27}$$

## 3.4 Geometric Interpretation of the Second-order Cone Model

An ellipsoid and interior of it centered at $\boldsymbol{\mu}_{\boldsymbol{\phi}}^{(i)}$, shape matrix $\Sigma_{\boldsymbol{\phi}}^{(i)}$ and radius $r$ is denoted by $\mathcal{E}(\boldsymbol{\mu}_{\boldsymbol{\phi}}^{(i)}, \Sigma_{\boldsymbol{\phi}}^{(i)}, r)$ as:

$$\mathbf{z} \in \mathcal{E}\left(\boldsymbol{\mu}_{\boldsymbol{\phi}}^{(i)}, \Sigma_{\boldsymbol{\phi}}^{(i)}, r\right) \triangleq \left\{\mathbf{z} \in \mathbb{R}^d : \left(\mathbf{z} - \boldsymbol{\mu}_{\boldsymbol{\phi}}^{(i)}\right)^T \left(\Sigma_{\boldsymbol{\phi}}^{(i)}\right)^{-1} \left(\mathbf{z} - \boldsymbol{\mu}_{\boldsymbol{\phi}}^{(i)}\right) \leq r^2\right\}$$

Equivalently we have:

$$\mathbf{z} \in \mathcal{E}\left(\boldsymbol{\mu}_{\boldsymbol{\phi}}^{(i)}, \Sigma_{\boldsymbol{\phi}}^{(i)}, r\right) \triangleq \left\{\mathbf{z} \in \mathbb{R}^d : \mathbf{z} = \boldsymbol{\mu}_{\boldsymbol{\phi}}^{(i)} + r\left(\Sigma_{\boldsymbol{\phi}}^{(i)}\right)^{\frac{1}{2}} \mathbf{u}, \|\mathbf{u}\|_2 \leq 1\right\}$$

Constraint (26) holds if and only if

$$\mathbf{w}^T \mathbf{z} \geq b - \xi_i \quad \forall \mathbf{z} \in \left\{\mathbf{z}^{(i)} \in \mathbb{R}^d : \mathbf{z}^{(i)} = \boldsymbol{\mu}^{(i)} + \sqrt{\frac{1-\alpha}{\alpha}} \left(\Sigma_{\boldsymbol{\phi}}^{(i)}\right)^{\frac{1}{2}} \mathbf{u}, \|\mathbf{u}\|_2 \leq 1\right\}$$

This shows that in the case of uncertainty, each data point $\tilde{\mathbf{x}}^{(i)}$ instead of a single point, takes its values from an ellipsoidal uncertainty set $\mathcal{E}(\boldsymbol{\mu}_{\boldsymbol{\phi}}^{(i)}, \Sigma_{\boldsymbol{\phi}}^{(i)}, r)$. This interprets model (17)-(19) as a robust optimization problem over an ellipsoidal uncertainty set. If every point $\mathbf{z}$ from the ellipsoid $\mathcal{E}(\boldsymbol{\mu}_{\boldsymbol{\phi}}^{(i)}, \Sigma_{\boldsymbol{\phi}}^{(i)}, r)$ satisfies the constraint $\mathbf{w}^T \mathbf{z} \geq b - \xi_i$, the second-order constraint (26) holds. As can be seen, the size of the uncertainty ellipsoid depends on the parameter $\alpha$. As $\alpha$ gets closer to zero, the size of the ellipsoid and consequently the uncertainty increases.

## 3.5 Dual of the Second-order Cone Model

By defining dual vector $\gamma$ and notation $\mathbf{z}(\mathbf{u}^{(i)}) = \boldsymbol{\mu}_{\boldsymbol{\phi}}^{(i)} + \sqrt{\frac{1-\alpha}{\alpha}} \left(\Sigma_{\boldsymbol{\phi}}^{(i)}\right)^{\frac{1}{2}} \mathbf{u}^{(i)}$ for $i = 1, \ldots, l$, the dual of model (25)-(27) is as follows:

$$\min_{\boldsymbol{\gamma} \in \mathbb{R}^l, \mathbf{u}^{(i)} \in \mathbb{R}^d} \frac{1}{2} \sum_{i=1}^{l} \sum_{i'=1}^{l} \gamma_i \gamma_{i'} \mathbf{z}(\mathbf{u}^{(i)})^T \mathbf{z}(\mathbf{u}^{(i')})$$

$$\sum_{i=1}^{l} \gamma_i = 1 \tag{28}$$

$$\|\mathbf{u}^{(i)}\|_2 \leq 1 \quad \forall i \in \{1, \ldots, l\}$$

$$0 \leq \gamma_i \leq \frac{1}{l\nu} \quad \forall i \in \{1, \ldots, l\}$$

Proof. The Lagrangian function of model (26) is:

$$\mathcal{L}(\mathbf{w}, b, \boldsymbol{\xi}, \boldsymbol{\gamma}, \boldsymbol{\beta}) = \frac{1}{2} \|\mathbf{w}\|_2^2 - b + \frac{1}{l\nu} \sum_{i=1}^{l} \xi_i$$
$$+ \sum_{i=1}^{l} \gamma_i \left( \sqrt{\frac{1-\alpha}{\alpha}} \left\| \left(\Sigma_{\boldsymbol{\phi}}^{(i)}\right)^{\frac{1}{2}} \mathbf{w} \right\|_2 - \mathbf{w}^T \boldsymbol{\mu}_{\boldsymbol{\phi}}^{(i)} - \xi_i + b \right) - \sum_{i=1}^{l} \beta_i \xi_i$$

According to the definition, of 2-norm, for any $\mathbf{x} \in \mathbb{R}^n$, $\|\mathbf{x}\|_2 = \max_{\|\mathbf{y}\|_2 \leq 1} \mathbf{y}^T \mathbf{x}$. Now we use this to eliminate the term $\left\| \left(\Sigma_{\boldsymbol{\phi}}^{(i)}\right)^{\frac{1}{2}} \mathbf{w} \right\|_2$ in the Lagrangian. By considering the vector $\mathbf{u}^{(i)}$ with $\|\mathbf{u}^{(i)}\|_2 \leq 1$, the Lagrangian is rewritten as (Here, the arbitrary vector $-\mathbf{u}^{(i)}$ has been used in the definition of 2-norm):

$$\mathcal{L}(\mathbf{w}, b, \boldsymbol{\xi}, \boldsymbol{\gamma}, \boldsymbol{\beta}) = \frac{1}{2} \|\mathbf{w}\|_2^2 - b + \frac{1}{l\nu} \sum_{i=1}^{l} \xi_i$$
$$+ \sum_{i=1}^{l} \gamma_i \left( \sqrt{\frac{1-\alpha}{\alpha}} \max_{\|-\mathbf{u}^{(i)}\|_2 \leq 1} -\mathbf{u}^{(i)T} \left(\Sigma_{\boldsymbol{\phi}}^{(i)}\right)^{\frac{1}{2}} \mathbf{w} - \mathbf{w}^T \boldsymbol{\mu}_{\boldsymbol{\phi}}^{(i)} - \xi_i + b \right) - \sum_{i=1}^{l} \beta_i \xi_i$$
$$= \max_{\|\mathbf{u}^{(i)}\|_2 \leq 1, \, i=1,\ldots,l} \frac{1}{2} \|\mathbf{w}\|_2^2 - b + \frac{1}{l\nu} \sum_{i=1}^{l} \xi_i$$
$$- \sum_{i=1}^{l} \gamma_i \left( \sqrt{\frac{1-\alpha}{\alpha}} (\mathbf{u}^{(i)})^T \left(\Sigma_{\boldsymbol{\phi}}^{(i)}\right)^{\frac{1}{2}} \mathbf{w} + \mathbf{w}^T \boldsymbol{\mu}_{\boldsymbol{\phi}}^{(i)} + \xi_i - b \right) - \sum_{i=1}^{l} \beta_i \xi_i$$

By defining:

$$\mathcal{L}_1(\mathbf{w}, b, \boldsymbol{\xi}, \boldsymbol{\gamma}, \boldsymbol{\beta}, \mathbf{u}^{(1)}, \ldots, \mathbf{u}^{(l)})$$
$$= \frac{1}{2}\|\mathbf{w}\|_2^2 - b + \frac{1}{l\nu}\sum_{i=1}^{l}\xi_i$$
$$- \sum_{i=1}^{l}\gamma_i\left(\sqrt{\frac{1-\alpha}{\alpha}}(\mathbf{u}^{(i)})^T\left(\Sigma_{\boldsymbol{\phi}}^{(i)}\right)^{\frac{1}{2}}\mathbf{w} + \mathbf{w}^T\boldsymbol{\mu}_{\boldsymbol{\phi}}^{(i)} + \xi_i - b\right) - \sum_{i=1}^{l}\beta_i\xi_i$$

the Lagrangian becomes:

$$\mathcal{L}(\mathbf{w}, b, \boldsymbol{\xi}, \boldsymbol{\gamma}, \boldsymbol{\beta}) = \max_{\|\mathbf{u}^{(i)}\|_2 \leq 1 \ i=1,\ldots,l} \mathcal{L}_1(\mathbf{w}, b, \boldsymbol{\xi}, \boldsymbol{\gamma}, \boldsymbol{\beta}, \mathbf{u}^{(1)}, \ldots, \mathbf{u}^{(l)})$$

Now the Lagrange dual function is defined as:

$$\mathcal{G}(\boldsymbol{\gamma}, \boldsymbol{\beta}) = \min_{\mathbf{w},b,\boldsymbol{\xi}} \mathcal{L}(\mathbf{w}, b, \boldsymbol{\xi}, \boldsymbol{\gamma}, \boldsymbol{\beta}) = \min_{\mathbf{w},b,\boldsymbol{\xi}} \max_{\|\mathbf{u}^{(i)}\|_2 \leq 1 \ i=1,\ldots,l} \mathcal{L}_1(\mathbf{w}, b, \boldsymbol{\xi}, \boldsymbol{\gamma}, \boldsymbol{\beta}, \mathbf{u}^{(1)}, \ldots, \mathbf{u}^{(l)})$$
$$= \max_{\|\mathbf{u}^i\|_2 \leq 1 \ i=1,\ldots,l} \min_{\mathbf{w},b,\boldsymbol{\xi}} \mathcal{L}_1(\mathbf{w}, b, \boldsymbol{\xi}, \boldsymbol{\gamma}, \boldsymbol{\beta}, \mathbf{u}^{(1)}, \ldots, \mathbf{u}^{(l)})$$

Since $\mathcal{L}_1$ is convex, by taking partial derivatives with respect to $\mathbf{w}$, $\boldsymbol{\xi}$, and b and solving the following system, we obtain its minimum value of $\mathcal{L}_1$:

$$\mathbf{w} = \sum_{i=1}^{l}\gamma_i \mathbf{z}(\mathbf{u}^{(i)})$$
$$\sum_{i=1}^{l}\gamma_i = 1 \tag{29}$$
$$\frac{1}{l\nu} - \gamma_i - \beta_i = 0$$

By substituting equations (29) in the Lagrangian function, we obtain the Lagrange dual function:

$$\mathcal{G}(\boldsymbol{\gamma}, \boldsymbol{\beta}) = \frac{1}{2}\|\mathbf{w}\|_2^2 - b + \frac{1}{l\nu}\sum_{i=1}^{l}\xi_i - \sum_{i=1}^{l}\gamma_i\left(\sqrt{\frac{1-\alpha}{\alpha}}(\mathbf{u}^{(i)})^T\left(\Sigma_{\boldsymbol{\phi}}^{(i)}\right)^{\frac{1}{2}}\mathbf{w} + \mathbf{w}^T\boldsymbol{\mu}_{\boldsymbol{\phi}}^{(i)} + \xi_i - b\right)$$
$$- \sum_{i=1}^{l}\beta_i\xi_i$$
$$= \frac{1}{2}\|\mathbf{w}\|_2^2 + \frac{1}{l\nu}\sum_{i=1}^{l}\xi_i - \sum_{i=1}^{l}\gamma_i\xi_i - \sum_{i=1}^{l}\beta_i\xi_i$$
$$- \sum_{i=1}^{l}\gamma_i\left(\sqrt{\frac{1-\alpha}{\alpha}}(\mathbf{u}^{(i)})^T\left(\Sigma_{\boldsymbol{\phi}}^{(i)}\right)^{\frac{1}{2}} + \left(\boldsymbol{\mu}_{\boldsymbol{\phi}}^{(i)}\right)^T\right)\mathbf{w} + \left(\sum_{i=1}^{l}\gamma_i - 1\right)b$$
$$= -\frac{1}{2}\sum_{i=1}^{l}\sum_{i'=1}^{l}\gamma_i\gamma_{i'}\mathbf{z}(\mathbf{u}^{(i)})^T\mathbf{z}(\mathbf{u}^{(i')})$$

Therefore, the dual objective function becomes:

$$\mathcal{G}(\boldsymbol{\gamma}, \boldsymbol{\beta}) = \max_{\|\mathbf{u}^i\|_2 \leq 1} \min_{\mathbf{w}, b, \boldsymbol{\xi}} \mathcal{L}_1(\mathbf{w}, b, \boldsymbol{\xi}, \boldsymbol{\gamma}, \boldsymbol{\beta}) = \max_{\|\mathbf{u}^i\|_2 \leq 1, \, i=1,\ldots,l} -\frac{1}{2} \sum_{i=1}^{l} \sum_{i'=1}^{l} \gamma_i \gamma_{i'} \mathbf{z}(\mathbf{u}^{(i)})^T \mathbf{z}(\mathbf{u}^{(i')})$$

Finally, according to KKT conditions, model (28) is the dual problem.

As can be seen, compared to the dual problem of the standard kernel-based one-class SVM (8), in the dual objective function (28), the nominal value $\mathbf{x}^{(i)}$ has been replaced with ellipsoid $\sqrt{\frac{1-\alpha}{\alpha}} \left(\Sigma_{\boldsymbol{\phi}}^{(i)}\right)^{\frac{1}{2}} \mathbf{u}^{(i)} + \boldsymbol{\mu}_{\boldsymbol{\phi}}^{(i)}$ where $\mathbf{u}^{(i)}$ is a vector with $\|\mathbf{u}^{(i)}\|_2 \leq 1$, which indicates the separation of ellipsoids $\mathcal{E}\left(\boldsymbol{\mu}_{\boldsymbol{\phi}}^{(i)}, \Sigma_{\boldsymbol{\phi}}^{(i)}, \sqrt{\frac{1-\alpha}{\alpha}}\right)$ for $i = 1, \ldots, l$ from the origin.

Suppose $\mathbf{w}^*, b^*, \boldsymbol{\xi}^*$ are the optimal solutions of primal model (25)-(27) and $\boldsymbol{\gamma}^*$ is the optimal value of dual vector, KKT conditions can be derived as:

Stationarity conditions:

$$\mathbf{w}^* = \sum_{i=1}^{l} \gamma_i^* \left( \sqrt{\frac{1-\alpha}{\alpha}} \left(\Sigma_{\boldsymbol{\phi}}^{(i)}\right)^{\frac{1}{2}} \left(\mathbf{u}^{(i)}\right)^* + \boldsymbol{\mu}_{\boldsymbol{\phi}}^{(i)} \right) \tag{30}$$

$$\sum_{i=1}^{l} \gamma_i^* = 1 \tag{31}$$

Primal Feasibility conditions:

$$\sqrt{\frac{1-\alpha}{\alpha}} \left\| \left(\Sigma_{\boldsymbol{\phi}}^{(i)}\right)^{\frac{1}{2}} \mathbf{w}^* \right\|_2 \leq \mathbf{w}^{*T} \boldsymbol{\mu}_{\boldsymbol{\phi}}^{(i)} + \xi_i^* - b \quad i = 1, \ldots, l \tag{32}$$

$$\xi_i^* \geq 0 \quad i = 1, \ldots, l$$

Dual Feasibility conditions:

$$0 \leq \gamma_i^* \leq \frac{1}{l\nu} \tag{33}$$

Complementary Slackness conditions:

$$\gamma_i^* \left( \sqrt{\frac{1-\alpha}{\alpha}} \left\| \left(\Sigma_{\boldsymbol{\phi}}^{(i)}\right)^{\frac{1}{2}} \mathbf{w}^* \right\|_2 - \mathbf{w}^{*T} \boldsymbol{\mu}_{\boldsymbol{\phi}}^{(i)} - \xi_i^* + b \right) = 0 \tag{34}$$

$$\left(\frac{1}{l\nu} - \gamma_i^*\right) \xi_i^* = 0$$

According to equations (30) and (31), the vector $\mathbf{w}^*$ is generated by a convex combination of the ellipsoids under uncertainty $\mathcal{E}\left(\boldsymbol{\mu}_{\boldsymbol{\phi}}^{(i)}, \Sigma_{\boldsymbol{\phi}}^{(i)}, \sqrt{\frac{1-\alpha}{\alpha}}\right)$ for $i = 1, \ldots, l$.

Similar to standard support vector machine models, we can define the support ellipsoid for model (28); $\mathcal{E}\left(\boldsymbol{\mu}_{\boldsymbol{\phi}}^{(i)}, \Sigma_{\boldsymbol{\phi}}^{(i)}, \sqrt{\frac{1-\alpha}{\alpha}}\right)$ is a support ellipsoid when $\gamma_i^* \neq 0$.

However, even with knowing $\mathbf{w}^*$, it is difficult to compute $\boldsymbol{\mu}_{\boldsymbol{\phi}}^i$ and $\Sigma_{\boldsymbol{\phi}}^i$ directly based on the two first moments $\boldsymbol{\mu}^i$ and $\Sigma^i$ in most cases, especially since the mapping $\boldsymbol{\phi}$ is usually defined implicitly and the mapped feature space can be infinite-dimensional $\mathbb{R}^\infty$.

## 3.6 Distributionally Robust Chance-Constrained Kernel-based One-Class SVM

As mentioned before, using the mapping function allows us to map a nonlinear dataset to a higher dimensional space such that the dataset becomes origin-inseparable in the new space. This mapping enables the model to represent more complex patterns. The main challenge is that after mapping the data to the feature space, exact distribution information about the mapped data points such as the mean and covariance is no longer available.

To address this issue, a data-driven approach has been proposed. In this approach, first, for each data point, some samples are generated from its distribution. These samples are mapped to the feature space through the mapping function. Using the mapped samples in feature space, the empirical mean and covariance matrices are estimated. Finally, using these empirical estimates, the chance constraints are approximated. However, the number of variables used in this reformulated model depends on the feature space dimension, and using a mapping that induces an infinite-dimensional space leads to an extremely large-sized problem. To overcome this issue, by adopting the kernel idea and using empirical estimates of the mean vector and covariance matrix, a kernel-based second-order cone model is provided as an approximation of model (25)-(27). Utilizing the kernel enables performing computations in the original space and makes the computational complexity independent of dimension of feature space. The research conducted in this subsection is inspired by [23].

Assume for each $i = 1, \ldots, l$, there are $m_i$ independent observations $S^{(i)} \triangleq \{\mathbf{x}^{(i_t)} \in \mathbb{R}^n : t = 1, \ldots, m_i\}$ for the uncertain input $\tilde{\mathbf{x}}^i$. Let $m = \sum_{i=1}^N m_i$ be the total number of samples $S = \{\mathbf{x}^{(s)}\}_{s=1}^m = \bigcup_{i=1}^l \{\mathbf{x}^{(i_t)}\}_{t=1}^{m_i}$. Also assume the random vectors $\tilde{\mathbf{x}}^i$ for $i = 1, \ldots, l$ are independent. Due to the independence assumption between $\tilde{\mathbf{x}}^i$s, we focus on one data $\tilde{\mathbf{x}}^i$ in the following. The main mechanism of data-driven approach to solve problem (27) is using $\boldsymbol{\phi}(S^{(i)}) \triangleq \{\boldsymbol{\phi}(\mathbf{x}^{(i_t)}) \in \mathbb{R}^d : t = 1, \ldots, m_i\}$ to estimate the moments $\boldsymbol{\mu}_{\boldsymbol{\phi}}^{(i)}$ and $\Sigma_{\boldsymbol{\phi}}^{(i)}$. For this purpose, the following approximations are used:

$$\hat{\boldsymbol{\mu}}_{\boldsymbol{\phi}}^{S^{(i)}} = \frac{1}{m_i} \sum_{t=1}^{m_i} \boldsymbol{\phi}(\mathbf{x}^{(i_t)}) \tag{35}$$

$$\widehat{\Sigma}_{\boldsymbol{\phi}}^{S^{(i)}} = \frac{1}{m_i}\sum_{t=1}^{m_i}\left(\boldsymbol{\phi}(\mathbf{x}^{(i_t)})-\widehat{\boldsymbol{\mu}}_{\boldsymbol{\phi}}^{S^{(i)}}\right)\left(\boldsymbol{\phi}(\mathbf{x}^{(i_t)})-\widehat{\boldsymbol{\mu}}_{\boldsymbol{\phi}}^{S^{(i)}}\right)^T$$
$$= \frac{1}{m_i}\sum_{t=1}^{m_i}\boldsymbol{\phi}(\mathbf{x}^{(i_t)})\boldsymbol{\phi}(\mathbf{x}^{(i_t)})^T - \widehat{\boldsymbol{\mu}}_{\boldsymbol{\phi}}^{S^{(i)}}\left(\widehat{\boldsymbol{\mu}}_{\boldsymbol{\phi}}^{S^{(i)}}\right)^T \quad (36)$$

According to equation (30) the vector $\mathbf{w}$ is a linear combination of points from the uncertainty ellipsoids $\mathcal{E}\left(\boldsymbol{\mu}_{\boldsymbol{\phi}}^{(i)}, \Sigma_{\boldsymbol{\phi}}^{(i)}, \sqrt{\frac{1-\alpha}{\alpha}}\right)$, $i \in \{1, \dots, l\}$. When the shape matrix of uncertainty ellipsoid $\mathcal{E}\left(\boldsymbol{\mu}_{\boldsymbol{\phi}}^{(i)}, \Sigma_{\boldsymbol{\phi}}^{(i)}, \sqrt{\frac{1-\alpha}{\alpha}}\right)$ corresponding to the $i^{th}$ data is determined by the estimated covariance matrix $\widehat{\Sigma}_{\boldsymbol{\phi}}^{(i)}$, any point in this ellipsoid is a linear combination of the mapped samples used in estimating the covariance matrix, since the eigenvectors of covariance matrix $\widehat{\Sigma}_{\boldsymbol{\phi}}^{(i)}$ lie in the space spanned by the mapped members of $S^{(i)}$ in feature space and the eigenvectors span the entire ellipsoid [19]. Hence, a member of estimated uncertainty ellipsoid $\mathcal{E}\left(\widehat{\boldsymbol{\mu}}_{\boldsymbol{\phi}}^{(i)}, \widehat{\Sigma}_{\boldsymbol{\phi}}^{S_i}, \sqrt{\frac{1-\alpha}{\alpha}}\right)$ is specified as:

$$\sqrt{\frac{1-\alpha}{\alpha}}\left(\widehat{\Sigma}_{\boldsymbol{\phi}}^{S_i\frac{1}{2}}\right)^T \mathbf{u}^i + \widehat{\boldsymbol{\mu}}_{\boldsymbol{\phi}}^{(i)} = \sum_{t=1}^{m_i} a_{i_j}^0 \boldsymbol{\phi}(\mathbf{x}^{i_t}) \quad (37)$$

By considering the following definitions:
$$\mathbf{p} = [p_{1_1} \quad \dots \quad p_{1_{m_1}} \quad \dots \quad p_{l_1} \quad \dots \quad p_{l_{m_l}}]^T \in \mathbb{R}^m \text{ where } p_{i_t} = \gamma_i a_{i_t}^0$$
$$\boldsymbol{\phi}(X) = [\boldsymbol{\phi}(\mathbf{x}^{(1_1)}) \quad \dots \quad \boldsymbol{\phi}(\mathbf{x}^{(1_{m_1})}) \quad \dots \quad \boldsymbol{\phi}(\mathbf{x}^{(l_1)}) \quad \dots \quad \boldsymbol{\phi}(\mathbf{x}^{(l_{m_l})})] \in \mathbb{R}^{d \times m}$$
the equation of $\mathbf{w}$ is rewritten as:

$$\mathbf{w} = \sum_{i=1}^{l}\gamma_i \sum_{t=1}^{m_i} a_{i_j}^0 \boldsymbol{\phi}(\mathbf{x}^{(i_t)}) = \boldsymbol{\phi}(X)\mathbf{p} \quad (38)$$

An important conclusion is the linear combination representation of the normal vector $\mathbf{w}$ in terms of the mapped samples that are involved in estimating the means and covariances.

By defining the kernel function $\begin{cases} k: \mathbb{R}^n \times \mathbb{R}^n \to \mathbb{R} \\ k(\mathbf{x}, \mathbf{x}') = \boldsymbol{\phi}(\mathbf{x})^T \boldsymbol{\phi}(\mathbf{x}') \end{cases}$ and the kernel matrix $\overline{K} \triangleq \boldsymbol{\phi}(X)^T \boldsymbol{\phi}(X)$ with elements

$$\overline{K}_{t,s} = \boldsymbol{\phi}(\mathbf{x}^{(t)})^T \boldsymbol{\phi}(\mathbf{x}^{(s)}) = k(\mathbf{x}^{(t)}, \mathbf{x}^{(s)}) \quad \forall t, s \in \{1, \dots, m\} \quad (39)$$

the following equations hold:

$$\mathbf{w}^T\mathbf{w} = (\boldsymbol{\phi}(X)\mathbf{p})^T\boldsymbol{\phi}(X)\mathbf{p} = \mathbf{p}^T\bar{K}\mathbf{p} \tag{40}$$

$$\mathbf{w}^T\hat{\boldsymbol{\mu}}_{\boldsymbol{\phi}}^{S^{(i)}} = (\boldsymbol{\phi}(X)\mathbf{p})^T\left(\frac{1}{m_i}\sum_{t=1}^{m_i}\boldsymbol{\phi}(\mathbf{x}^{(i_t)})\right) = \mathbf{p}^T\left(\frac{1}{m_i}\boldsymbol{\phi}(X)^T\sum_{t=1}^{m_i}\boldsymbol{\phi}(\mathbf{x}^{(i_t)})\right) = \mathbf{p}^T\bar{\mathbf{k}}^{(i)} \tag{41}$$

where

$$\bar{\mathbf{k}}^{(i)} = \frac{1}{m_i}\left[\sum_{t=1}^{m_i}k(\mathbf{x}^{(1_1)},\mathbf{x}^{(i_t)}) \quad \ldots \quad \sum_{t=1}^{m_i}k(\mathbf{x}^{(l_{m_l})},\mathbf{x}^{(i_t)})\right]_{m\times 1}^T \tag{42}$$

$$\mathbf{w}^T\hat{\Sigma}_{\boldsymbol{\phi}}^{S_i}\mathbf{w} = \mathbf{w}^T\left(\frac{1}{m_i}\sum_{t=1}^{m_i}\left(\boldsymbol{\phi}(\mathbf{x}^{(i_t)})-\hat{\boldsymbol{\mu}}_{\boldsymbol{\phi}}^{S^{(i)}}\right)\left(\boldsymbol{\phi}(\mathbf{x}^{(i_t)})-\hat{\boldsymbol{\mu}}_{\boldsymbol{\phi}}^{S^{(i)}}\right)^T\right)\mathbf{w} = \mathbf{p}^T\bar{K}^{(i_t)}\mathbf{p} \tag{43}$$

where $\bar{K}^{(i_t)} = \left[k(\mathbf{x}^{(1_1)},\mathbf{x}^{(i_t)}) \quad \ldots \quad k(\mathbf{x}^{(l_{m_l})},\mathbf{x}^{(i_t)})\right]_{m\times 1}^T$.

By defining the matrix

$$\Sigma_K^{(i)} = \frac{1}{m_i}\sum_{t=1}^{m_i}(\bar{K}^{(i_t)}-\bar{\mathbf{k}}^{(i)})(\bar{K}^{(i_t)}-\bar{\mathbf{k}}^{(i)})^T \in \mathbb{R}^{m\times m} \tag{44}$$

and substituting equations (40), (41) and (43) in model (25)-(27), a second-order cone model based on kernel is obtained as an approximation of model (25)-(27):

$$\begin{aligned}\min_{\mathbf{p},b\geq 0,\boldsymbol{\xi}} \quad & \frac{1}{2}\mathbf{p}^T\bar{K}\mathbf{p} - b + \frac{1}{l\nu}\sum_{i=1}^{l}\xi_i \\ & \sqrt{\frac{1-\alpha}{\alpha}}\left\|\left(\Sigma_K^{(i)}\right)^{\frac{1}{2}}\mathbf{p}\right\|_2 \leq \mathbf{p}^T\bar{\mathbf{k}}^{(i)} + \xi_i - b \quad i=1,\ldots,l \\ & \xi_i \geq 0 \quad i=1,\ldots,l\end{aligned} \tag{45}$$

After finding the optimal solution $\mathbf{p}^*$, $b^*$, and $\boldsymbol{\xi}^*$ of model (40), the decision function is determined as:

$$f_{\boldsymbol{\phi}}(\mathbf{x};\mathbf{v}^*,b^*) = -sign\left(\sum_{i=1}^{l}\sum_{j=1}^{m_i}\mathbf{p}_{i_j}^*k(\mathbf{x}^{(i_j)},\mathbf{x}) - b^*\right) \tag{46}$$

The number of constraints in model (45) equals the number of training data points and the size of square matrix $\left(\Sigma_K^{(i)}\right)^{\frac{1}{2}}$ equals the total number of samples. We call the approach of using the entire original dataset and generating samples for each data point $KDRCC - Sampling$. To reduce model size and the size of matrix $\left(\Sigma_K^{(i)}\right)^{\frac{1}{2}}$, instead of generating new samples, the dataset can be clustered using a clustering method and each cluster can be regarded as a sample. More precisely, after determining the clusters, the closest data point to the center of each cluster is considered as the mean corresponding to that cluster and members of that are considered as the generated samples corresponding to that cluster. This approach is called $KDRCC - Clustering\ I$. After determining the clusters, cluster members other than the mean can be removed and then by an estimated covariance matrix, for example the covariance matrix estimated using cluster members, samples corresponding to that cluster are generated. This approach is called $KDRCC - Clustering\ II$. In both last approaches, the number of constraints equals the number of clusters which is less than the total number of data points. Also, in the second approach, the size of matrix $\left(\Sigma_K^{(i)}\right)^{\frac{1}{2}}$ equals the total number of data points which is less than the number of generated samples in the first approach. Moreover, by properly selecting the number of samples per cluster in the third approach, the size of matrix $\left(\Sigma_K^{(i)}\right)^{\frac{1}{2}}$ will be less conservative compared to the corresponding one in the first approach.

## 4 Evaluation Metrics and Computational Experiments

This section contains the computational experiments conducted on the proposed model (45). All computational experiments are implemented using the CVXPY library in Python and solved by Mosek on a desktop equipped with Intel(R) Core(TM) i5-4200U CPUs @ 2.6GHz and 16 GB RAM.

The evaluation metrics including accuracy (Acc), Positive F1 score (PF1) as the harmonic mean of precision and recall for the anomaly class, and Negative F1 (NF1) score for the normal class are used to evaluate model performance.

It is assumed that the first and second-order moment information is available for each random vector $\tilde{\mathbf{x}}^i$. To evaluate the robustness of the proposed model (45), its performance is evaluated on data sets generated by distributions such as normal, uniform, and Student's t (with degree of freedom 7) with the same moments as the given means and covariances in terms of various evaluation metrics. Also, a separate validation set is used to determine hyperparameters and the overall performance of the model on a separate test set from the training and validation sets is evaluated and compared to the deterministic one-class SVM model (8).

The following steps are taken to apply model (45) on the datasets:

The training data points and their means and covariances are given.

For the $i^{th}$ training data, $n_{batch}$ samples are generated following an arbitrary distribution with its corresponding mean and covariance matrix.

Based on the generated samples, the kernel matrix $\bar{K}$ with elements (39), vector $\bar{\mathbf{k}}^{(i)}$ (42) and matrix $\Sigma_K^{(i)}$ (44) are formed.

Model (45) is formulated using the matrices $\bar{K}$, and $\Sigma_K^{(i)}$, and the vector $\bar{\mathbf{k}}^{(i)}$.

The hyperparameters $\nu$ and $\gamma$ are selected using the validation set by the holdout selection method and positive F1-score metric [25].

After finding the optimal solution of model (45), the decision function is determined by equation (46).

Different evaluation metrics are calculated based on the true and predicted labels on the set set.

First, a 2D dataset ($D_1$) is utilized to gain geometric intuition about the decision boundary obtained from model (45). The normal data points are generated from a normal distribution with mean $[2,2]^T$ and covariance matrix $0.3 I_n$ where $I_n$ is the n-dimensional identity matrix. The anomalies are generated with mean $[0.2,0.2]^T$ and covariance $0.1 I_n$. The training set includes 100 normal points and 2 anomalies as outliers. The validation and test sets both contain 35 normal points and 5 anomalies each, randomly selected. It is assumed the first and second-order moment information is available. Specifically, the available data points in the training set are considered as the mean and for simplicity the covariance matrix is taken as $0.1 I_n$.

Figure 3 shows the Training partition of dataset $D_1$. The blue points represent normal data and the red ones are anomalies.

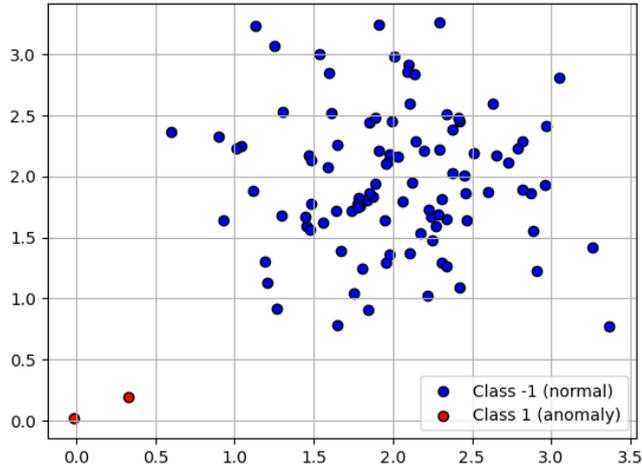

Figure 3: Training partition of dataset $D_1$

The results obtained with $\alpha = 0.01$, $n_{batch} = 5$, $\nu = 0.37$ and $\gamma = 0.00098$ on the test partition of dataset $D_1$ is reported in Table 1:

Table 1: Performance of the proposed model (45) with $KDRCC - sampling$ approach and (8) on the test set $D_1$

| Evaluation metric | Normal distribution | Uniform distribution | Student's t distribution |
|---|---|---|---|
| Acc (%) | 99.65 | 99.90 | 99.12 |
| PF1 (%) | 90.90 | 97.43 | 77.5 |
| NF1 (%) | 99.82 | 99.95 | 99.55 |

As can be seen, the decision boundary obtained from model (45) can maintain proper values with respect to different evaluation metrics across the three considered distributions. This validates that model (45) can be robust against distribution uncertainty.

The mean points (solid color circles) and samples generated under normal, uniform and t distributions (light color circles) as well as the obtained optimal decision boundary from model (45) are shown in Figure 4:

Figure 4: Test datasets with normal, uniform and Student's t distributions along with the optimal decision boundary of model (45)

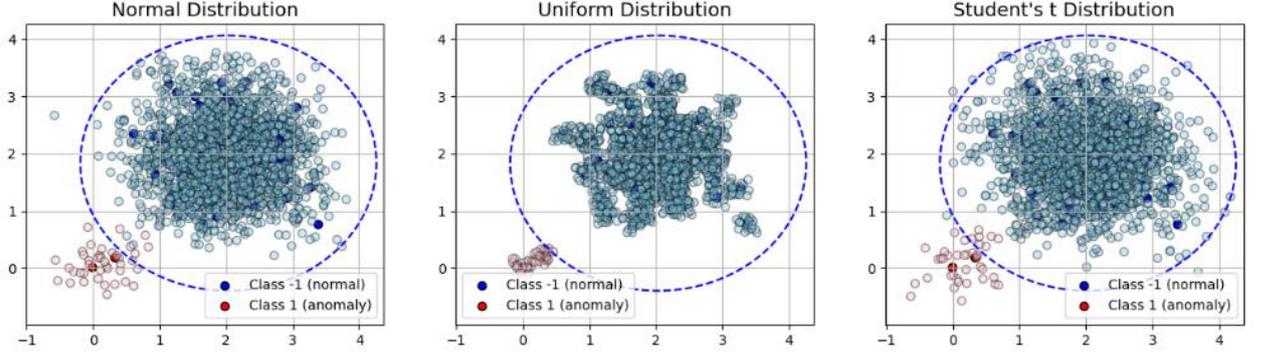

Now we compare the results obtained from the standard SVM model (8) and the proposed model (45). In both models, the holdout method and the positive F1-score metric are used to select hyperparameters. In addition to dataset $D_1$, the model is tested on two other simulated datasets. Dataset $D_2$ has 7 features where the features 1, 5, and 6 are class-discriminative features while the features 2, 3, 4, and 7 are not. Feature 1 has a normal distribution with mean -6 and variance 1 in the normal class and has a normal distribution with mean 6 and variance 1 in the anomaly class. Feature 2 has a normal distribution with mean -3 and variance 5 in both classes. Feature 3 is a quasi-constant feature with mean 10 and variance 0.001 in both classes. Feature 4 has normal distributions with means 1 and 1.02 and variance 0.008 in the normal and anomaly classes respectively. Feature 5 is twice feature 1. Feature 6 has uniform distributions over [-5, 5] and [30, 50] in the normal and anomaly classes respectively. Feature 7 has a uniform distribution over [-50, 50] in both classes. The class ratios in the training, validation, and test sets of $D_2$ are similar to the corresponding sets in $D_1$. The training data points are considered as the means and the Identity matrix is taken as the covariance matrix. Dataset $D_3$ includes normal data points sampled from the parametric curves $x_2 = r\cos(t)$ and $x_1 = r\sin(t)$ where $0 \leq t \leq 2\pi$ and $3 \leq r \leq 5$ and anomalies sampled from the same parametric equations with parameter $10 \leq r \leq 12$. The class ratios and chosen means and covariance matrix are similar to $D_2$.

The results of the proposed models (45) with three approaches $KDRCC - Sampling$, $KDRCC - Clustering$ I and $KDRCC - Clustering$ II and the standard one-class SVM model (8) on the test set in terms of accuracy (Acc), positive F1 score (PF1), negative F1 score (NF1), and F1 score on the sets generated by three distributions (NUT-F1) are reported in Table 2 and Table 3. Also, for more reliable results, each experiment is done with changing data partitioning into training, validation and test sets by altering the seed parameter in random data selection. boldfaced numbers indicate the best values obtained by a model on a dataset based on an evaluation metric.

Table 2: Performance of the proposed model (45) with $KDRCC-sampling$ approach and standard one-class SVM model (8)

| Dataset | $KDRCC-Sampling$ (45) | | | | One-Class SVM model (8) | | | |
|---|---|---|---|---|---|---|---|---|
| | Acc | PF1 | NF1 | NUT-PF1 | Acc | PF1 | NF1 | NUT-PF1 |
| $D_1$ | **100** | **100** | **100** | N: **90.90**<br>U: **97.43**<br>T: 77.5 | 97.5 | 90.90 | 98.55 | N: 19.85<br>U: 30.89<br>T: 17.32 |
| | **100** | **100** | **100** | N: 89.19<br>U: **100**<br>T: 79.01 | 97.5 | 90.90 | 98.55 | N: 22.66<br>U: 32.79<br>T: 18.60 |
| | **100** | **100** | **100** | N: 90.90<br>U: 100<br>T: 76.54 | 97.5 | 90.90 | 98.55 | N: 44.31<br>U: 68.96<br>T: 39.40 |
| $D_2$ | **100** | **100** | **100** | N: 89.89<br>U: **100**<br>T: 88.89 | **100** | **100** | **100** | N: 35.24<br>U: 34.48<br>T: 35.71 |
| | **100** | **100** | **100** | N: **100**<br>U: **100**<br>T: **100** | 85 | 62.5 | 90.62 | N: 26.75<br>U: 27.59<br>T: 27.40 |

| | | | | N:<br>**100** | | | | N:<br>17.35 |
| --- | --- | --- | --- | --- | --- | --- | --- | --- |
| | **100** | **100** | **100** | U:<br>**100** | 85 | 62.5 | 90.62 | U:<br>21.0 |
| | | | | T:<br>**100** | | | | T:<br>16.56 |
| | | | | N:<br>**96.10** | | | | N:<br>7.61 |
| | **100** | **100** | **100** | U:<br>**100** | 70 | 45.45 | 79.31 | U:<br>8.96 |
| | | | | T:<br>**92.10** | | | | T:<br>7.82 |
| | | | | N:<br>97.43 | | | | N:<br>9.10 |
| $D_3$ | **100** | **100** | **100** | U:<br>**100** | 75 | 50 | 83.33 | U:<br>12.10 |
| | | | | T:<br>**97.5** | | | | T:<br>9.03 |
| | | | | N:<br>**100** | | | | N:<br>9.65 |
| | **100** | **100** | **100** | U:<br>**100** | 75 | 50 | 83.33 | U:<br>12.65 |
| | | | | T:<br>**97.5** | | | | T:<br>9.65 |

Table 3: Performance of the proposed model (45) with $KDRCC - Clustering\ I$ and $KDRCC - Clustering\ II$ approaches

| Dataset | $KDRCC - Clustering\ I$ | | | | $KDRCC - Clustering\ II$ | | | |
| --- | --- | --- | --- | --- | --- | --- | --- | --- |
| | Acc | PF1 | NF1 | NUT-PF1 | Acc | PF1 | NF1 | NUT-PF1 |
| $D_1$ | **100** | **100** | **100** | N:<br>85.05<br><br>U:<br>88.89<br><br>T:<br>70.70 | **100** | **100** | **100** | N:<br>90.66<br><br>U:<br>94.74<br><br>T:<br>**86.11** |
| | **100** | **100** | **100** | N: **92.10** U: **100** T: **79.52** | **100** | **100** | **100** | N: 76.92 U: 91.89 T: 70.27 |
| | **100** | **100** | **100** | N: 82.97 U: 100 T: 75.51 | **100** | **100** | **100** | N: 83.54 U: 100 T: 83.54 |
| $D_2$ | **100** | **100** | **100** | N: 90.90 U: **100** T: 86.95 | **100** | **100** | **100** | N: **100** U: **100** T: **94.12** |
| | **100** | **100** | **100** | N: **100** U: **100** T: **100** | **100** | **100** | **100** | N: **100** U: **100** T: **100** |
| | **100** | **100** | **100** | N: **100** U: **100** T: 97.56 | 100 | 100 | 100 | N: **100** U: **100** T: 98.7**6** |
| $D_3$ | **100** | **100** | **100** | N: 94.87 U: **100** T: 83.72 | **100** | **100** | **100** | N: **96.10** U: **100** T: 90.48 |

|   |     |     |     | N: 100        |     |     |     | N: 97.43     |
|---|-----|-----|-----|---------------|-----|-----|-----|--------------|
|   | 100 | 100 | 100 | U: 100        | 100 | 100 | 100 | U: 100       |
|   |     |     |     | T: 91.76      |     |     |     | T: **97.5**  |
|   |     |     |     | N: 100        |     |     |     | N: **100**   |
|   | 100 | 100 | 100 | U: 100        | 100 | 100 | 100 | U: **100**   |
|   |     |     |     | T: 86.96      |     |     |     | T: **100**   |

By comparing the results of the proposed model (45) and the standard model (8), it is observed that for all test sets and sets generated by the three distributions, model (45) achieves better values in terms of various evaluation metrics. The superiority in robustness against uncertainty in different distributions is more significant.

## 5 Conclusion and Future Work

In this paper, a nonlinear one-class SVM model is presented by considering uncertainty in the feature vector of each data point and its probability distribution. This model is proposed to guarantee desirable performance in the worst case over different distributions, enable classification of non-separable datasets, and reduce computational complexity through the kernel idea. The computational results demonstrate the superiority of the proposed model compared to the standard one-class SVM in terms of evaluation metrics on various test sets. It is also observed that the proposed model has a relatively consistent and satisfactory performance across different distributions indicating its robustness against uncertainty in data distributions. For future work, it is suggested to study applying chance constraints simultaneously for all data points through joint chance constraints and providing an efficient solution method for large-scale problems.